\documentstyle[12pt]{amsart}

\newtheorem{theo+}           {Theorem}      [section]
\newtheorem{prop+}  [theo+]  {Proposition}
\newtheorem{coro+}  [theo+]  {Corollary}
\newtheorem{lemm+}  [theo+]  {Lemma}
\newtheorem{exam+}  [theo+]  {Example}
\newtheorem{rema+}  [theo+]  {Remark}
\newtheorem{defi+}  [theo+]  {Definition}

\newenvironment{theorem}{\begin{theo+}}{\end{theo+}}
\newenvironment{proposition}{\begin{prop+}}{\end{prop+}}
\newenvironment{corollary}{\begin{coro+}}{\end{coro+}}
\newenvironment{lemma}{\begin{lemm+}}{\end{lemm+}}
\newenvironment{example}{\begin{exam+}\normalshape}{\end{exam+}}
\newenvironment{remark}{\begin{rema+}\normalshape}{\end{rema+}}

\def \rn{\Bbb R}
\def \cn{\Bbb C}

\def \H{\cal H}

\def \V{\cal V}

\def \sp#1#2{\langle #1,#2 \rangle}

\def \nab#1#2{\hbox{$\nabla$\kern -.3em\lower 1.0 ex
     \hbox{$#1$}\kern -.1 em {$#2$}}}
\def \nabm#1#2{\hbox{$\nabla$\kern -.3em\lower 1.0 ex
     \hbox{$#1$}\kern -.1 em {$#2$}}}
\def \nabn#1#2{\hbox{$\nabla^N$\kern -.3em\lower 1.0 ex
     \hbox{$#1$}\kern -.1 em {$#2$}}}
\def \nabh#1#2{\hbox{$\nabla^{\H}$\kern -.3em\lower 1.0 ex
     \hbox{$#1$}\kern -.1 em {$#2$}}}
\def \nabl#1#2{\hbox{$\nabla^L$\kern -.3em\lower 1.0 ex
     \hbox{$#1$}\kern -.1 em {$#2$}}}
\def \nabf#1#2{\hbox{$\nabla^{\phi^{-1}TN}$\kern -.3em\lower 1.0 ex
     \hbox{$#1$}\kern -.1 em {$#2$}}}
\def \nabfg#1#2{\hbox{$\nabla^{\phi^{-1}TG_n({\cn}^{m+n})}$
\kern -.3em\lower 1.0 ex \hbox{$#1$}\kern -.1 em {$#2$}}}

\def \proc#1{{\cn} P^{#1}}

\def \SU#1{\text{\bf SU}(#1)}
\def \U#1{\text{\bf U}(#1)}

\begin{document}
\baselineskip 20pt

\title{Harmonic Morphisms between \\
       Almost Hermitian Manifolds}

\author{SIGMUNDUR GUDMUNDSSON}
\author{JOHN C. WOOD}

\thanks{The first author was supported by
        The Swedish Natural Science Research Council.}
\thanks{Both authors were partially supported by EC grant
CHRX-CT92-0050.}

\date{Dedicato alla memoria del professore Franco Tricerri}

\keywords{almost Hermitian manifolds,
          harmonic maps,
          harmonic morphisms}

\subjclass{58E20}

\address{Department of Mathematics,
         University of Lund,
         Box 118,
         S-221 00 Lund,
         Sweden}

\email{sigma@@maths.lth.se}

\address{Department of Pure Mathematics,
         University of Leeds,
         Leeds LS2 9JT,
         UK}

\email{j.c.wood@@leeds.ac.uk}

\maketitle

\def\abstractname{Sunto}%
\begin{abstract}
In questo articolo si ottengono condizioni sulla forma di Lee affich\`e
un' applicazione olomorfa tra variet\`a quasi-hermitiane sia un'
applicazione armonica oppure un morfismo armonico.  Poi, si discute
sotto quali condizioni, (i) l'immagine di un' applicazione olomorfa di
una variet\`a cosimplettica \`e anche cosimplettica, (ii) un'
applicazione olomorfa con immagine hermitiana definisce una struttura
hermitiana sul suo dominio.
\end{abstract}

\section{Introduction}

One of Franco Tricerri's interests was in Hermitian manifolds.
In \cite{T-V} various types of Hermitian structures are
discussed and  conditions on the Lee form are of paramount importance.
That Hermitian structures are closely connected with harmonic
morphisms is shown in \cite{B-W} and \cite{W-4d}.  In this paper we
study this connection for more general almost Hermitian manifolds.
We obtain conditions involving the Lee form under which holomorphic
maps between almost Hermitian manifolds are harmonic maps or
morphisms.  We show that the image of certain holomorphic maps from a
cosymplectic manifold is cosymplectic if and only if the map is a
harmonic morphism, generalizing a result of Watson \cite{Wat}.
Finally, in Theorem \ref{theo-integral} we give conditions under which
a harmonic morphism into a Hermitian manifold defines an integrable
Hermitian structure on its domain.

\section{Harmonic morphisms}

For a smooth map $\phi:(M,g)\to(N,h)$ between Riemannian manifolds
its {\it tension field} $\tau(\phi)$ is the trace of the second
fundamental form $\nabla d\phi$ of $\phi$:
\begin{equation}\label{equa-tau}
\tau(\phi)=\sum_{j}\{\nabf{e_j}{d\phi(e_j)}-d\phi(\nabm{e_j}{e_j})\}
\end{equation}
where $\{e_j\}$ is a local orthonormal frame for $TM$,
$\nabla^{\phi^{-1}TN}$
denotes the pull-back of the Levi-Civita connection $\nabla^N$ on $N$ to the
pull-back bundle  $\phi^{-1}TN\to M$ and $d\phi:TM\to\phi^{-1}TN$
is the pull-back of the differential of $\phi$.
The map $\phi$ is said to be {\it harmonic} if its
tension field vanishes i.e. $\tau(\phi)=0$.  J. Eells and J. H. Sampson
proved in \cite{E-S} that any holomorphic map between
K\"ahler manifolds is harmonic and this was later generalized by
A. Lichnerowicz in \cite{L}.  For information on harmonic maps, see
\cite{E-L-1}, \cite{E-L-2}, \cite{E-L-3} and the references therein.

A {\it harmonic morphism} is a smooth map $\phi:(M,g)\to(N,h)$
between Riemannian manifolds which pulls back germs of
real-valued harmonic functions on $N$ to germs of harmonic functions
on $M$.  A smooth map $\phi:M\to N$ is called {\it horizontally
(weakly) conformal} if for each $x\in M$ {\it either}
\begin{enumerate}
  \item[(i)] the rank of the differential $d\phi_x$ is $0$, {\it or}
  \item[(ii)] for the orthogonal decomposition
$T_xM={\H}_x\oplus{\V}_x$ with ${\V}_x=\ker\ {d\phi}_x$ the restriction
$d\phi_x|_{{\H}_x}$ is a conformal linear map {\it onto} $T_{\phi(x)}N$.
\end{enumerate}
Points of type (i) are called {\it critical points} of $\phi$ and those of
type (ii) {\it regular points}.
The conformal factor $\lambda(x)$ is called the {\it dilation} of
$\phi$ at $x$.  Setting $\lambda=0$ at the critical points gives a
continuous function $\lambda:M\to [0,\infty)$ which is smooth at regular
points, but whose  square $\lambda^2$ is smooth on the whole of $M$.
Note that at a regular point $\phi$ is a submersion.
A horizontally weakly conformal map is called {\it horizontally homothetic}
if $d\phi(\text{grad}(\lambda^2))=0$.  B. Fuglede showed in \cite{F-2} that
a horizontally homothetic harmonic morphism has no critical points.

The following characterization of harmonic morphisms is due to
Fuglede and T. Ishihara, see \cite{F-1}, \cite{I}:  {\it A smooth
map $\phi$ is a harmonic morphism if and only if it is a
horizontally weakly conformal harmonic map}.  More geometrically we
have the following result due to P. Baird and Eells, see \cite{B-E}:

\begin{theorem}\label{theo-B-E}
Let $\phi:(M,g)\to(N,h)$ be a non-constant horizontally weakly
conformal map.  Then
  \begin{enumerate}
    \item[(i)] if $N$ is a surface, i.e. of real dimension $2$, then $\phi$
      is a harmonic morphism
      if and only if its fibres are minimal at regular points,
    \item[(ii)] if the real dimension of $N$ is greater than $2$
         then any two of the following conditions imply the third:
         \begin{enumerate}
           \item[(a)] $\phi$ is a harmonic morphism,
           \item[(b)] the fibres of $\phi$ are minimal at regular points,
           \item[(c)] $\phi$ is horizontally homothetic.
         \end{enumerate}
  \end{enumerate}
\end{theorem}

Harmonic morphisms exhibit many properties which are ``dual'' to those
of harmonic maps.  For example, whereas harmonic maps exhibit conformal
invariance in a $2$-dimensional domain (cf. \cite{E-S}, Proposition
p.126), harmonic morphisms exhibit conformal invariance in a
$2$-dimensional codomain:  {\it If $\phi:(M,g)\to(N,h)$ is a
harmonic morphism to a $2$-dimensional Riemannian manifold and
$\psi:(N,h)\to(\tilde N,\tilde h)$ is a weakly conformal map to
another $2$-dimensional Riemannian manifold, then the composition
$\psi\circ\phi$ is a harmonic morphism}.  In particular the concept of
a {\it harmonic morphism to a Riemann surface} is well-defined.  For
information on harmonic morphisms see \cite{B-W} and \cite{W-Sendai}.

\section{Almost Hermitian manifolds}

Let $(M^m,g,J)$ be an almost Hermitian manifold i.e. a Riemannian
manifold $(M,g)$ of even real dimension $2m$ together with
an almost complex structure $J:TM\to TM$ which is isometric
on each tangent space and satisfies $J^2=-I$.  Let $T^{\cn}M$ be the
complexification of the tangent bundle.  We then have an
orthogonal decomposition $$T^{\cn}M=T^{1,0}M\oplus T^{0,1}M$$ of
$T^{\cn}M$ into the $\pm i$-eigenspaces of $J$, respectively.
Each vector $X\in T^{\cn}M$ can be written as $X=X^{1,0}+X^{0,1}$ with
$$X^{1,0}=\frac 12(X-iJX)\in T^{1,0}M\ \ \text{and}\ \
  X^{0,1}=\frac 12(X+iJX)\in T^{0,1}M,$$
and locally one can always choose an orthonormal frame
$\{e_1,\dots,e_m,Je_1,\dots,Je_m\}$ for $TM$ such that
$$T^{1,0}M=\text{span}_{\cn}
\{Z_1=\frac{e_1-iJe_1}{\sqrt 2},\dots,
  Z_m=\frac{e_m-iJe_m}{\sqrt 2}\},$$
$$T^{0,1}M=\text{span}_{\cn}
\{\bar Z_1=\frac{e_1+iJe_1}{\sqrt 2},\dots,
  \bar Z_m=\frac{e_m+iJe_m}{\sqrt 2}\}.$$
The set $\{Z_k|\ k=1,\dots,m\}$ is called a local {\it Hermitian frame}
on $M$.

As for any other $(1,1)$-tensor the {\it divergence} of $J$ is given by
\begin{eqnarray*}\delta J=\text{div}(J)
&=&\sum_{k=1}^{m}(\nabm{e_k}J)(e_k)+(\nabm{Je_k}J)(Je_k)\\
&=&\sum_{k=1}^{m}(\nabm{\bar Z_k}J)(Z_k)+(\nabm{Z_k}J)(\bar Z_k).
\end{eqnarray*}

\begin{remark}  Modulo a constant, the vector field $J\delta J$ is
  the dual to the Lee form, see \cite{T-V}.  It is called the
  {\it Lee vector field}.
\end{remark}

Following Kot\=o \cite{K} and Gray \cite{Gr-2} with alternative
terminology due to Salamon \cite{S} we call an almost Hermitian manifold
$(M,g,J)$
\begin{enumerate}
  \item[(i)] {\it quasi-K\"ahler} or {\it $(1,2)$-symplectic} if
$$(\nabm XJ)Y+(\nabm{JX}J)JY=0\ \ \text{for all $X,Y\in C^\infty(TM)$},\
\ \text{and}$$
  \item[(ii)] {\it semi-K\"ahler} or {\it cosymplectic} if $\delta J=0$.
\end{enumerate}
Note that a $(1,2)$-symplectic manifold $(M,g,J)$ is automatically
cosymplectic.  It is an easy exercise to prove the following two
well-known results:

\begin{lemma}\label{lemm-(1,2)-symplectic}
  Let $(M,g,J)$ be an almost Hermitian manifold.  Then the following
  conditions are equivalent:
  \begin{enumerate}
\item[(i)] $M$ is $(1,2)$-symplectic,
\item[(ii)] $\nab{\bar Z}W\in C^\infty(T^{1,0}M)$ for all
  $Z,W\in C^\infty(T^{1,0}M)$.
  \end{enumerate}
\end{lemma}

\begin{lemma}\label{lemm-cosymplectic}
Let $(M,g,J)$ be an almost Hermitian manifold.  Then the following
conditions are equivalent:
\begin{enumerate}
\item[(i)] $M$ is cosymplectic,
\item[(ii)] $\sum_{k=1}^m\nab{\bar Z_k}{Z_k}\in C^\infty(T^{1,0}M)$
for any local Hermitian frame $\{Z_k\}$.
\end{enumerate}
\end{lemma}

\begin{example}\label{exam-Cal-Eck}
For $r,s\ge 0$ let $(M,g)$ be the product
$S^{2r+1}\times S^{2s+1}$ of the two unit spheres in ${\cn}^{r+1}$
and ${\cn}^{s+1}$ equipped with their standard Euclidean metrics.  The
manifold $(M,g)$ has a standard almost Hermitian structure $J$ which can
be described as follows (cf. \cite{Gr-Her} and \cite{Tri-Van}):
Let $n_1,n_2$ be the
unit normals to $S^{2r+1},S^{2s+1}$ in ${\cn}^{r+1},{\cn}^{s+1}$ and let
${\H}_1,{\H}_2$ be the horizontal spaces of the Hopf maps
$S^{2r+1}\to\proc r$,
$S^{2s+1}\to\proc s$, respectively.  Then any vector tangent to
$M$ has the form $$X=X_1+aJ_1n_1+X_2+bJ_2n_2$$ where $a,b\in\rn$,
$X_1\in{\H}_1$, $X_2\in{\H}_2$, and $J_1,J_2$ are the standard K\"ahler
structures on ${\cn}^{r+1}$ and ${\cn}^{s+1}$, respectively.
Then the almost complex structure $J$
on $M$ is given by $$J:X\mapsto J_1X_1-bJ_1n_1+J_2X_2+aJ_2n_2.$$
We calculate that (cf. \cite{CMW}):
$$\delta J=-2(rJ_1n_1+sJ_2n_2).$$
The almost Hermitian manifold $(M,g,J)$ is called the {\it Calabi-Eckmann
manifold}.  It is cosymplectic if and only if $s=r=0$ i.e. when
$M$ is the real $2$-dimensional torus in ${\cn}^2$.
\end{example}

\begin{example}\label{exam-twistor}
Any invariant metric on a $3$-symmetric space gives it a
$(1,2)$-symplec\-tic structure (cf. Proposition 3.2 of \cite{Gr-3}).
Such $3$-symmetric spaces occur as twistor spaces of symmetric spaces.
One interesting example is the complex
Grassmannian $G_n({\cn}^{m+n})=\SU{m+n}/\text{\bf S}(\U m\times \U n)$
with twistor bundle the flag manifold
$N=\SU{m+n}/\text{\bf S}(\U m\times \U k\times\U{n-k})$ and
projection $\pi:N\to G_n({\cn}^{m+n})$ induced by the inclusion map
$\U k\times\U{n-k}\hookrightarrow\U n$.  The manifold $N$ has an almost
Hermitian structure usually denoted by $J^2$ such that
$(N,g,J^2)$ is $(1,2)$-symplectic for any $\SU{m+n}$-invariant metric
$g$ on $N$.  For further details see \cite{S}.
\end{example}

Finally, recall that an almost Hermitian manifold is called {\it
Hermitian} if its almost complex structure is integrable.
A necessary and sufficient condition for this is the vanishing of
the Nijenhuis tensor (cf. \cite{N-N}), or equivalently, that
$T^{1,0}M$ is closed under the Lie bracket i.e.
$[T^{1,0}M,T^{1,0}M]\subset T^{1,0}M$.

\section{The harmonicity of holomorphic maps}

Throughout this section we shall assume that $(M^m,g,J)$ and
$(N^n,h,J^N)$ are almost Hermitian manifolds of complex dimensions
$m$ and $n$ with Levi-Civita connections $\nabla$ and $\nabla^N$,
respectively.  Furthermore we suppose that the map $\phi:M\to N$ is
holomorphic i.e. its differential $d\phi$ satisfies
$d\phi\circ J=J^N\circ d\phi$.  We are interested in studying under
what additional assumptions the map $\phi$ is a harmonic map or morphism.

For a local Hermitian frame $\{Z_k\}$ on $M$ we define
$$A=\sum_{k=1}^m\nabf{\bar Z_k}d\phi(Z_k)\ \ \text{and}\ \
  B=-\sum_{k=1}^md\phi(\nabm{\bar Z_k}{Z_k}).$$

\begin{lemma}\label{lemm-tension}
  Let $\phi:(M,g,J)\to(N,h,J^N)$ be a
  map between almost Hermitian manifolds.  If $N$ is $(1,2)$-symplectic
  then the tension field $\tau(\phi)$ of $\phi$ is given by
  \begin{equation}\label{equa-tau-J}
  \tau(\phi)=-d\phi(J\delta J).
  \end{equation}
\end{lemma}

\begin{proof}
Let $\{Z_k\}$ be a local Hermitian frame, then a simple calculation
shows that
$$J\delta J
=\sum_{k=1}^m\{(1+iJ)\nabm{\bar Z_k}{Z_k}+(1-iJ)\nabm{Z_k}{\bar Z_k}\},$$
so that the $(0,1)$-part of $J\delta J$ is given by
$$(J\delta J)^{0,1}=(1+iJ)\sum_{k=1}^m\nabm{\bar Z_k}{Z_k}.$$
The holomorphy of $\phi$ implies that $d\phi(Z_k)$ belongs
to $C^\infty(\phi^{-1}T^{1,0}N)$ and the $(1,2)$-symplecticity
on $N$ that $\nabf{\bar Z_k}d\phi(Z_k)\in C^\infty(\phi^{-1}T^{1,0}N)$.
This means that $A^{0,1}=0$.
{}From Equation (\ref{equa-tau}) and the symmetry of the second
fundamental form $\nabla d\phi$ we deduce that $\tau(\phi)=2(A+B)$.
Taking the $(0,1)$-part and using the
fact that $\phi$ is holomorphic we obtain
$$\tau(\phi)^{0,1}=2(A^{0,1}+B^{0,1})=2B^{0,1}=-d\phi(J\delta J)^{0,1}.$$
  Since $\tau(\phi)$ and $d\phi(J\delta J)$ are both real, we deduce
  the result.
\end{proof}
\vskip .5cm

The next proposition gives a criterion for harmonicity in terms of the Lee
vector field.

\begin{proposition}\label{prop-harmonic}
Let $\phi:(M,g,J)\to(N,h,J^N)$ be a holomorphic map from an
almost Hermitian manifold to a $(1,2)$-symplectic manifold.
Then $\phi$ is harmonic if and only if $d\phi(J\delta J)=0$.
\end{proposition}

Note that since we are assuming that $\phi$ is holomorphic
$d\phi(J\delta J)=0$ is equivalent to $d\phi(\delta J)=0$.
As a direct consequence of Proposition \ref{prop-harmonic}
we have the following result of Lichnerowicz, see \cite{L}:

\begin{corollary}\label{corr-L}
  Let $\phi:(M,g,J)\to(N,h,J^N)$ be a holomorphic map from a
  cosymplectic manifold to a $(1,2)$-symplectic one.  Then $\phi$
  is harmonic.
\end{corollary}

To deduce that $\phi$ is a harmonic morphism we must assume that $\phi$
is horizontally weakly conformal.  In that situation we can say more:

\begin{proposition}\label{prop-harm-morph-1}
  Let $\phi:(M^{m},g,J)\to(N^{n},h,J^N)$ be a surjective horizontally weakly
  conformal holomorphic map between almost Hermitian manifolds.
  Then any two of the following conditions imply the third:
  \begin{enumerate}
    \item[(i)] $\phi$ is harmonic and so a harmonic morphism,
    \item[(ii)] $d\phi(J\delta J)=0$.
    \item[(iii)] $N$ is cosymplectic,

  \end{enumerate}
\end{proposition}

\begin{proof}  By taking the $(0,1)$-part of equation (\ref{equa-tau})
we obtain  $$\tau(\phi)^{0,1}=2(A^{0,1}+B^{0,1}).$$
The tension field $\tau(\phi)$ is real so that the map
$\phi$ is harmonic if and only if $\tau(\phi)^{0,1}=0$.  Since
$2B^{0,1}=-d\phi(J\delta J)^{0,1}$ and the vector field
$d\phi(J\delta J)$ is real the condition
$d\phi(J\delta J)=0$ is equivalent to $B^{0,1}=0$.  To complete the proof
we shall now show that $A^{0,1}=0$ on $M$ if and only if $N$ is
cosymplectic.

Let $R$ be the open subset of regular points of $\phi$.  Let
$p\in R$ and $\{Z'_1,\dots,Z'_n\}$ a local Hermitian frame
on an open neighbourhood $V$ of $\phi(p)\in N$.  Let
$Z^*_1,\dots,Z^*_n$ be the unique horizontal lifts of $Z'_1,\dots,Z'_n$ to
$\phi^{-1}(V)$ and normalize by setting
$Z_k=\lambda Z^*_k$ for $k=1,2,\dots,n$, where $\lambda$ is the
dilation of $\phi$ defined in Section 2.  Then we can, on an open
neighbourhood of $p$, extend $\{Z_1,\dots,Z_n\}$ to a local Hermitian frame
$\{Z_1,\dots,Z_m\}$ for $M$.  We then have
$$A=\sum_{k=1}^n\nabf{\bar Z_k}{(\lambda Z'_k)}=
\sum_{k=1}^n\bar Z_k(\lambda)Z'_k
+\lambda^2\sum_{k=1}^n\nabn{\bar Z'_k}{Z'_k}.$$
The vector field $\sum_{k=1}^n\bar Z_k(\lambda)Z'_k$ belongs to
$\phi^{-1}T^{1,0}N$, so by Lemma \ref{lemm-cosymplectic}, $A^{0,1}$
vanishes on $R$ if and only if $N$ is cosymplectic at each point of
$\phi(R)$.

Now note that if $p$ is a critical point of $\phi$ then {\it either} $p$
is a limit point of a sequence of regular points {\it or} $p$ is
contained in an open subset $W$ of critical points.
In the first case, if $A^{0,1}$ vanishes on $R$ then it vanishes
also at $p$ by
continuity.  In the second case $d\phi=0$ on $W$ so that $A^{0,1}=0$ at $p$.
This means that $A^{0,1}$ vanishes on $R$ if and only if it vanishes on $M$.

On the other hand, since $\phi$ is surjective it follows from Sard's
theorem that $\phi(R)$ is dense in $N$.  This implies that
$N$ is cosymplectic at points of $\phi(R)$ if and only if $N$ is
cosymplectic everywhere.
Putting the above remarks together yields the proof.
\end{proof}
\vskip .5cm

As a direct consequence of Proposition \ref{prop-harm-morph-1} we
have the following:

\begin{theorem}\label{theo-harm-morp}
  Let $\phi:(M,g,J)\to (N,h,J^N)$ be a surjective horizontally weakly
  conformal holomorphic map from a cosymplectic manifold to an almost
  Hermitian manifold. Then $N$ is cosymplectic if and only if $\phi$
  is a harmonic morphism.
\end{theorem}

Combining Theorems \ref{theo-harm-morp} and \ref{theo-B-E} we then obtain:

\begin{corollary}\label{corr-harm-morp-2}
Let $\phi:(M,g,J)\to(N,h,J^N)$ be a surjective horizontally homothetic
holomorphic map from a cosymplectic manifold to an almost Hermitian
manifold.  Then $N$ is cosymplectic if and only if $\phi$ has minimal
fibres.
\end{corollary}

Corollary \ref{corr-harm-morp-2} generalizes a result of B. Watson
in \cite{Wat} where it is assumed that the map $\phi$ is a Riemannian
submersion.  If the manifold $(M,g,J)$ is $(1,2)$-symplectic we have the
following version of Theorem \ref{theo-harm-morp}:

\begin{proposition}\label{prop-harm-morph-2}
Let $\phi:(M,g,J)\to (N^n,h,J^N)$ be a horizontally weakly
conformal holomorphic map from a $(1,2)$-symplectic manifold
to a cosymplectic one.  Then $\phi$  is a harmonic
morphism whose fibres are minimal at regular points.
If $n>1$ then $\phi$ is horizontally homothetic.
\end{proposition}

\begin{proof}
The inclusion maps of the fibres of $\phi$ are holomorphic maps between
$(1,2)$-symplectic manifolds.  They are isometric immersions
and, by Corollary \ref{corr-L}, harmonic so the fibres
are minimal.   For an alternative argument see \cite{Gr-1}.  If $n>1$ then
Theorem \ref{theo-B-E} implies that $\phi$ is horizontally homothetic.
\end{proof}
\vskip .5cm

Now assume that $n=1$, then $N$ is automatically K\"ahler and therefore
$(1,2)$-symplectic.  Further, any holomorphic map from an almost Hermitian
manifold $(M,g,J)$ to $N$ is horizontally weakly conformal.
Hence Proposition \ref{prop-harmonic} implies the following results:

\begin{corollary}\label{corr-surf-1}
  Let $\phi:(M,g,J)\to N$ be a holomorphic map from an almost
  Hermitian manifold to a Riemann surface.  Then $\phi$ is a harmonic
  morphism if and only if $d\phi(J\delta J)=0$.
\end{corollary}

\begin{corollary}\label{corr-surf-2}
  Let $\phi:(M,g,J)\to N$ be a holomorphic map from a cosymplectic
  manifold to a Riemann surface. Then $\phi$ is a harmonic morphism.
\end{corollary}

\begin{example}
For two integers $r,s\ge 0$ let $M$ be the Calabi-Eckmann manifold
$(S^{2r+1}\times S^{2s+1},g,J)$ and $\phi:M\to\proc r\times\proc s$
be the product of the Hopf maps $S^{2r+1}\to\proc r$, $S^{2s+1}\to\proc s$.
Then it is not difficult to see that $\phi$ is holomorphic.
Further the
kernel of $d\phi$ is given by
$\ker d\phi=\text{span}\{J_1n_1,J_2n_2\}$.
{}From Example \ref{exam-Cal-Eck} we get
$d\phi(\delta J)=-2d\phi(rJ_1n_1+sJ_2n_2)=0$.
Since the map $\phi$ is a Riemannian submersion we deduce by Proposition
\ref{prop-harm-morph-1} that $\phi$ is a harmonic morphism.
\end{example}

The next result can be extended to any of the twistor spaces considered
by Salamon in \cite{S}, but for clarity we state it for a particular case.

\begin{proposition}\label{prop-twistor}
Let $\pi:N\to G_n({\cn}^{m+n})$ be the twistor
fibration of Example \ref{exam-twistor} and $\phi:(M,g,J)\to N$
be a holomorphic map from an almost Hermitian manifold into the flag
manifold $N$.  Although $\psi=\pi\circ\phi:M\to G_n({\cn}^{m+n})$ is not,
in general, a holomorphic map, we have $$\tau(\psi)=-d\psi(J\delta J).$$
\end{proposition}

\begin{proof}
Let $\{Z_k\}$ be a local Hermitian frame on $M$.  Then by using
the Composition Law for the tension field and
Lemmas \ref{lemm-tension} and \ref{lemm-vanish} we obtain:
\begin{eqnarray*}\tau(\psi)&=&d\pi(\tau(\phi))
+\sum_{k=1}^m\nabla d\pi(d\phi(\bar Z_k),d\phi(Z_k))\\
&=&-d\pi(d\phi(J\delta J))+0\\
&=&-d\psi(J\delta J).
\end{eqnarray*}
\end{proof}

\begin{lemma}\label{lemm-vanish}
The twistor fibration $\pi:N\to G_n({\cn}^{m+n})$ is $(1,1)$-geodesic
i.e. $$\nabla d\pi(Z,W)=0$$ for all $p\in N$, $Z\in T^{1,0}_pN$ and
$W\in T^{0,1}_pN$.
\end{lemma}

\begin{proof}
Decompose $Z$ and $W$ into vertical and horizontal parts
$Z=Z^{\V}+Z^{\H}$, $W=W^{\V}+W^{\H}$.  Now since $\pi$ is a Riemannian
submersion $\nabla d\pi(Z^{\H},W^{\H})=0$ by Lemma 1.3 of \cite{O}.  Further
$$\nabla d\pi(Z^{\V},W)=\nabfg{Z^{\V}}{d\pi(W)}
-d\pi(\nabm{Z^{\V}}W).$$
The first term is zero and the second term is of type $(0,1)$
with respect to the almost Hermitian structure $J_p$ on
$G_n({\cn}^{m+n})$ defined by $p$, since $\nabm{Z^{\V}}W$ is of type
$(0,1)$ and $d\pi_p:T_pM\to T_{\pi(p)}G_n({\cn}^{m+n})$ intertwines
$J$ and $J_p$.

Similarily $\nabla d\pi(W,Z^{\V})$ is of type $(1,0)$, so by the symmetry
of $\nabla d\pi$, $\nabla d\pi(Z^{\V},W)=0$.  Hence
$$\nabla d\pi(Z,W)=\nabla d\pi(Z^{\V},W)+\nabla d\pi(Z^{\H},W^{\H})
+\nabla d\pi(Z^{\H},W^{\V})=0.$$
\end{proof}

\section{Superminimality}

Let $\phi:(M,g)\to(N,h,J^N)$ be a horizontally conformal submersion
from a Riemannian manifold to an almost Hermitian manifold.
Assume that the fibres of $\phi$ are orientable and of real dimension $2$.
Then we can construct an almost Hermitian structure $J$ on $(M,g)$
such that $\phi$ becomes holomorphic:  make a smooth choice of an almost
Hermitian structure on each fibre and lift $J^N$ to the horizontal
spaces $\H$ using $d\phi\circ J=J^N\circ d\phi$.

For an almost Hermitian manifold $(M,g,J)$ we shall call an almost complex
submanifold $F$ of $M$ {\it superminimal} if $J$ is parallel
along $F$ i.e. $\nabm VJ=0$ for all vector fields $V$ tangent to $F$.
It is not difficult to see that any superminimal $F$ is minimal.
Superminimality of surfaces in $4$-dimensional manifolds has been
discussed by several authors, see for example \cite{Br}.

\begin{theorem}\label{theo-integral}
Let $\phi:(M,g,J)\to(N,h,J^N)$ be a horizontally conformal holomorphic
map from an almost Hermitian manifold to a Hermitian manifold with
complex $1$-dimensional fibres.  If
\begin{enumerate}
\item[(i)] the fibres of $\phi$ superminimal with respect to $J$, and
\item[(ii)] the horizontal distribution $\H$ satisfies
$[{\H}^{1,0},{\H}^{1,0}]^{\V}\subset {\V}^{1,0}$,
\end{enumerate}
then $J$ is integrable.
\end{theorem}

\begin{proof}
We will show that $T^{1,0}M$ is closed under the Lie bracket i.e.
$[T^{1,0}M,T^{1,0}M]\subset T^{1,0}M$, or equivalently:
\begin{enumerate}
\item[(a)] $[{\V}^{1,0},{\V}^{1,0}]\subset T^{1,0}M$,
\item[(b)] $[{\H}^{1,0},{\H}^{1,0}]\subset T^{1,0}M$,
\item[(c)] $[{\H}^{1,0},{\V}^{1,0}]\subset T^{1,0}M$.
\end{enumerate}
The fibres are complex $1$-dimensional so $[{\V}^{1,0},{\V}^{1,0}]=0$.
This proves (a).  Let $Z,W$ be two vector fields on $N$ of
type $(1,0)$ and let $Z^*,W^*$ be their horizontal lifts to ${\H}^{1,0}$.
Then $d\phi[Z^*,W^*]=[Z,W]$ is of type $(1,0)$ since $J^N$ is integrable.
The holomorphy of $\phi$ and assumption (ii) then imply (b).

Let $\sp{}{}$ be the complex bilinear extention of $g$ to $T^{\cn}M$ and
$V$ be a vertical vector field of type $(1,0)$.
Then $d\phi([V,Z^*])=[d\phi(V),Z]=0$, so $[V,Z^*]^{\H}=0$.
On the other hand,
\begin{eqnarray*}
\sp{[V,Z^*]}V&=&\sp{\nabm V{Z^*}}V-\sp{\nabm {Z^*}V}V\\
&=&-\sp{Z^*}{\nabm VV}-\frac 12Z^*(\sp VV).
\end{eqnarray*}
The subspace $T^{1,0}M$ is isotropic w.r.t. $\sp{}{}$ so $\sp VV=0$.
The superminimality of the fibres implies that
$J(\nabm VV)=\nabm V{JV}=i\nabm VV$.  Hence $\nabm VV$ is an
element of $T^{1,0}M$ so $\sp{Z^*}{\nabm VV}=0$.  This shows that
$\sp{[V,Z^*]}V=0$ so $[V,Z^*]^{\V}$ belongs to
$T^{1,0}M$.  This completes the proof.
\end{proof}
\vskip .5cm

The reader should note that condition (ii) of Theorem \ref{theo-integral}
is satisfied when the horizontal distribution $\H$ is integrable,
or when $N$ is complex $1$-dimensional, since in both cases
$[{\H}^{1,0},{\H}^{1,0}]^{\V}=0$.  Another example where the Theorem
\ref{theo-integral} applies is the following:

\begin{example}
The Hopf map $\phi:{\cn}^{n+1}-\{0\}\to\proc n$ is a horizontally
conformal submersion with complex $1$-dimensional fibres.
The horizontal distribution is non-integrable, but it is easily
seen that condition (ii) is satisfied, in fact
$[{\H}^{1,0},{\H}^{1,0}]^{\V}=0$.  The K\"ahler structure on $\proc n$
lifts to two almost Hermitian structures on ${\cn}^{n+1}-\{0\}$, the
fibres are superminimal with respect to both of these, so by Theorem
\ref{theo-integral} they are both Hermitian.  In fact one is the standard
K\"ahler structure; the other is not K\"ahler.
\end{example}

If $N$ is complex $1$-dimensional and $M$ is real $4$-dimensional
then Theorem \ref{theo-integral} reduces to a result of the second
author given in Proposition 3.9 of \cite{W-4d}.

\begin{example}
Let $(M^2,g,J)$ be an almost Hermitian manifold of complex dimension $2$
and $N$ a Riemann surface.  Then we have the identities \cite{Gau}
$$\nabm{\delta J}J=\nabm{J\delta J}J=0.$$  Otherwise said,
span$_{\rn}\{\delta J,J\delta J\}$ is contained in $\ker\nabla J$.  The
condition $d\phi(\delta J)=0$ for a non-constant holomorphic map
$\phi:M\to N$ is thus equivalent to the superminimality of the
fibres (at regular points) so that Corollary \ref{corr-surf-1}
translates into the following result of the second author,
see Proposition 1.3 of \cite{W-4d}.
{\it A holomorphic map from a Hermitian manifold of complex dimension
$2$ to a Riemann surface is a harmonic morphism if and only if its
fibres are superminimal at the regular points of $\phi$}.
\end{example}

\vskip 1cm

\end{document}